\documentclass[12pt]{article}
\usepackage[margin=1in]{geometry}
\usepackage{pdfpages}
\usepackage[utf8]{inputenc}
\usepackage{graphicx}
\usepackage{booktabs}
\usepackage{lipsum}
\usepackage{enumerate}
\usepackage{comment}
\usepackage{amssymb,amsthm,amsmath,bm}
\usepackage[colorlinks=true, urlcolor=blue,linkcolor=blue, citecolor=blue]{hyperref}

\title{Hedging Portfolio for a Degenerate Market Model\footnote{This work is supported by Tubitak Project No. 118F403.}}

\author{
	Mine \c{C}a\u{g}lar\\
	\texttt{Koç University}
	\and
 \.Ihsan Demirel\\
	\texttt{Koç University}
	\and
	Ali S\"{u}leyman \"{U}st\"{u}nel\\
   \texttt{Bilkent University}
}
\date{}

\newcommand{\nnumber}{\stepcounter{equation}
\tag{\theequation}\label}
	\newtheorem{theorem}{Theorem}[section]

	\newtheorem{proposition}{Proposition}[section]
	\newtheorem{lemma}{Lemma}[section]

	\newtheorem{remark}{Remark}[section]

\newcommand{\E}{\mathbb{E}}
\newcommand{\Q}{\mathbb{Q}}
\newcommand{\R}{\mathbb{R}}
\newcommand{\M}{\mathbb{M}}

\newcommand{\FF}{\mathcal{F}}

\begin{document}

\maketitle

\abstract{We consider a semimartingale market model when the underlying diffusion has a singular volatility matrix and compute the hedging portfolio for a given payoff function. Recently, the representation problem for such  degenerate diffusions with respect to a minimal martingale has been completely settled. This martingale representation and  Malliavin calculus established further for the functionals of a degenerate diffusion process  constitute the basis of the present work. Using the Clark-Hausmann-Bismut-Ocone type representation formula derived for these functionals, we prove a version of this formula under an equivalent martingale measure. This allows us to derive the hedging portfolio  as a  solution of a system of linear equations. The uniqueness of the solution is achieved by a projection idea that lies at the core of the martingale representation at the first place. We   demonstrate the hedging strategy as explicitly as possible with some examples of the payoff function such as those used in exotic options, whose value at maturity depends on the prices over the entire time horizon.  }


\section{Introduction}

An important application of the classical martingale representation theorem is in finance for calculating the hedging strategy when the risky asset price can be modeled as a diffusion process  with a strictly positive volatility. For a portfolio of assets that are diffusions in $\mathbb{R}^n$, the volatility is captured by a diffusion matrix $\sigma \in \mathbb{R}^{n\times d}$,  and the hedging  strategy is derived under the assumption that the  matrix $\sigma\sigma^*$ is non-singular.  The diffusion process is said to be non-degenerate in this case as studied extensively in prior work for hedging (see e.g. \cite{Klebaner,shreve2003}). On the other hand, recent advances in the martingale representation of degenerate diffusions \cite{ustunel1, ustunel2} make  the calculation of a hedging strategy possible when the volatility matrix $\sigma\sigma^*$ is singular in view of analysis on Wiener space \cite{ustunel_book_1995}. We show that the replicating portfolio process can be characterized as a solution of a system of linear equations based on these results.

Let $(\Omega, {\cal H}, \mathbb{P})$ be a probability space, and let $X=\{ X_{t}: 0\le t \le 1 \}$ satisfy the stochastic differential equation
\begin{equation} \label{diffusion}
    d X_{t}=b\left( X_t \right) d t + \sigma\left(X_t \right) d W_{t},
\end{equation}
where $\left\{W_{t}:0\le t\le 1\right\}$ is an $\mathbb{R}^{d}$-valued Brownian motion and $\sigma: \mathbb{R}^{n} \rightarrow \mathbb{R}^{n\times d}$ and $b: \mathbb{R}^{n} \rightarrow \mathbb{R}^{n}$ are measurable maps. We assume that the drift $b$ and the diffusion matrix $\sigma$ are Lipschitz and of linear growth as sufficient for $X$ to be the unique strong solution of \eqref{diffusion}. The diffusion $X$ is possibly degenerate, that is, $\sigma(x)\sigma(x)^*$ can be singular for some $x\in \R^n$.
Let $\mathcal{F}(X)= \{\mathcal{F}_t(X): 0\le t \le 1\}$ denote the filtration generated by $X$.  The martingale representation theorem  \cite[Thm.2]{ustunel1} reveals, in particular, that an $\mathcal{F}_1$-measurable functional $F$ of $X$ can be represented as
\[
F(X)=\mathbb{E}[F(X)] + \int_0^1 P(X_s) \xi_s(X) \cdot dW_s =
\mathbb{E}[F(X)] + \int_0^1 \xi_s(X) \cdot P(X_s) dW_s
\]
with  an $\mathcal{F}_t(X)$-adapted process  $\xi$ taking values in $\mathbb{R}^d$, where dot product is used for simplicity of notation and $ P(X_s)$ denotes orthogonal projection to the range space of $\sigma^*$, the transpose of $\sigma$. In essence, there exists a minimal martingale, given above as
$P(X_s) dW_s$ in its infinitesimal It\^o form, with respect to which every square integrable $\mathcal{F}_1$-measurable functional can be written as an integral of an $\mathcal{F}(X)$-adapted process. The representation problem for degenerate diffusions has been settled in \cite{ustunel1} as a result.

Further  in \cite{ustunel2},  Malliavin calculus for degenerate diffusions is developed, which  forms the basis for the results of the present paper. Let $(W, H, \mu)$ be the classical Wiener space on $\mathbb{R}^{d}$. For suitable $\mathcal{F}$-measurable functionals $F$, \cite[Thm.6]{ustunel2} provides a Clark-Hausmann-Bismut-Ocone type formula as
\begin{equation}  \label{Clark}
F(X)= \E[F(X)]+\int_0^1 P(X_s) \E[\hat{D}_s F(X)|{\cal F}_s] \cdot dW_s
\end{equation}
where the operator $\hat{D}$ is defined as the density of $\hat{\nabla}$ with respect to Lebesgue measure and $\hat{\nabla}$ is an operator analogous to  Gross-Sobolev derivative $\nabla$ for Wiener functionals. Starting with this formula, we consider the hedging of a stock portfolio when the prices are modeled as degenerate diffusions in this paper.

We not only solve the hedging problem for a semimartingale market model, but also find a  hedging strategy   to $\mathcal{F}(X)$, the filtration of the asset prices themselves, instead of the filtration $\mathcal{F}(W)$ of the driving Wiener process. More explicitly, let the price dynamics of $n$ assets $X_t$ in a market  follow  \eqref{diffusion} and let the equation for the risk-free asset $X^0_t$ at time $t$ be given by  $d X^{0}_t=r_t X^{0}_t d t$, where $r$ is the interest rate process, for $t\in [0,1]$.
 Let $\theta_t$ and $\theta^0_t$ be the number of shares of $n$ risky assets and the risk-free asset,  respectively,  where $\theta_t$ is taken as a row vector. Then, the value process $V^\theta_t$ is written as $
V^{\theta}_t= \theta_t\,X_t  + \theta^{0}_t\, X^{0}_t
$. Assuming that the portfolio $(\theta,\theta^0)$ is self-financing, we find the hedging  portfolio that replicates the terminal value function   $V_1$, which is assumed to be specified by an $\mathcal{F}$-measurable random variable $G(X)$, where $G$ is the payoff function. We take  $F=e^{-\int_0^1 r_s ds}G(X)$ in  \eqref{Clark} for our derivations of a hedging portfolio.
When there is no arbitrage, by denoting the equivalent martingale measure with $\Q$ and the $d$-dimensional Brownian motion under $\Q$ with $\widetilde{W}$, we show that the hedging strategy is obtained, in particular when $r$ is deterministic, by solving the equation
\begin{align*}
   \lefteqn{ \sigma^*(X_t) \theta_t^* =} \\
   & \quad \quad \displaystyle{ e^{-\int_t^1 r_s ds} P(X_t)\E_{\Q}\left[\hat{D}_{t} G(X) -G(X) \int_{t}^{1} P(X_s) \hat{D}_{t} \left(P(X_s)u(X_s)\right)\cdot d \widetilde{W}_s \mid \mathcal{F}_{t}(X)\right]}
\end{align*}
for $\theta$, which may not be unique although $\{P(X_t)\theta_t: 0\le t\le 1\}$ is unique for all solutions $\theta$.

We prove two fundamental results related to $\hat{\nabla}$  needed in our derivations. Namely, Proposition \ref{chain rule} as the chain rule, and Lemma \ref{LemmaFTC} as the fundamental theorem of calculus are developed as a follow up of \cite{ustunel2}, where $\hat{\nabla}$ is shown to satisfy the properties of a derivative operator.
Clearly, these properties are adopted from those of $\nabla$, but with care on the projection with $P(X_s)$ and using the cylindrical functions common in the domains of the two operators when necessary. In Theorem \ref{thm2}, we derive an equivalent representation to \eqref{Clark} using the equivalent martingale measure $\Q$ and  Wiener process $\widetilde{W}$, in view of the properties of the operator $\hat{\nabla}$.

As for applications, we limit ourselves to demonstrating the hedging strategy as explicitly as possible with some specific examples of the payoff $G$ as the discounted terminal value. Exotic options are considered in some detail as their value at maturity depends on the prices over the whole time horizon where Gross-Sobolev derivative is applicable.  The results can be useful in several finance and interdisciplinary applications where diffusion processes and hedging are considered (see e.g. \cite{ozekici2020,1999,Gobet2005,Hillairet2018,Tsoularis2018}).
 On the other hand, degeneracy in stochastic volatility models is investigated from the aspect of partial integro-differential equations that appear in corresponding risk-neutral pricing problems \cite{Costantini2012,Donatucci2016}. The volatility inherits randomness only from the prices $X$ and not stochastic on its own in the present paper, where the projection  $P(X_t)$ plays a crucial role in handling the degeneracy.

The paper is organized as follows. In Section 2, we review the essential parts of
Malliavin  calculus for degenerate diffusions and prove the preliminary results useful for the present work. Then, the hedging formula is derived for the degenerate semimartingale market model in Section 3. Special cases of the payoff function  are considered in Section 4 to demonstrate  hedging and option pricing.  Finally,  Section 5 concludes the paper.

\section{Preliminaries}

Let $\mathcal{S}(X)$ denote the set of cylindrical functions on the Wiener space $W$, given by
$$
    \mathcal{S}(X)=\left\{f\left(X_{t_{1}}^1, \ldots, X_{t_1}^n, \ldots, X_{t_m}^1, \ldots, X_{t_{m}}^n\right): 0 \leq t_{1}<\ldots<t_{m}, f \in \mathcal{S}\left(\mathbb{R}^{nm}\right), m \geq 1\right\}\
$$
where $\mathcal{S}\left(\mathbb{R}^{n}\right)$ denotes the space of rapidly decreasing smooth functions of Laurent Schwartz. In \cite{ustunel2}, for $h \in$ $H$ and $F(X) \in \mathcal{S}(X),$
an operator $\hat{\nabla}_h$  is defined  as
\begin{align*}
    \hat{\nabla}_h F(X)= \sum_{j=1}^m\sum_{i=1}^{n} \partial_{(j-1)n+i} f(X_{t_{1}}^1, \ldots, X_{t_1}^n, \ldots, X_{t_m}^1, \ldots, X_{t_{m}}^n)E[\nabla_{h}X_{t_j}^i|\mathcal{F}_1(X)].
\end{align*}
 where $\nabla$ denotes Gross-Sobolev derivative defined on Wiener space $(W,H,\mu).$ It has been shown in \cite[Cor.2]{ustunel2} that $\hat{\nabla}_{h}$ is a closable operator on  $L^{2}(\nu)$, where $\nu$ denotes the probability law of $X$, that is,  the image of $\mu$ under $X$.
 Then, also the operator $\hat{\nabla}$ can  be defined    by
\begin{align*}
    \hat{\nabla} F(X)=\sum_{i=1}^{\infty}  \hat{\nabla}_{e_i} F(X) e_i
\end{align*}
for $F(X) \in \mathcal{S}(X),$ where $\{e_i,\ i\geq 1\}$ is is an orthonormal basis in the Cameron-Martin space $H.$  It follows that $\hat{\nabla}$ is a closable operator from $L^{p}(\nu)$ to  $L^{p}(\nu;H)$, where $H$ is indicated to specify the range  of  $\hat{\nabla}$.

 The norm
$$
\|  F(X)\|_{p,1}:=  \|F(X)\|_{L^p(\mu)} + \| \hat{\nabla} F(X) \|_{L^p(\mu; H)}
$$
is used for the completion of $S(X)$, which will be denoted by $\M_{p,1}$. Note that  we use $|\cdot|$ for Euclidean norm,   $\| \cdot \|$ for $L^2([0,1])$-norm, and for all  others we specify the space in the notation.
For $F(X) \in \M_{2,1},$ define
 $\hat{D}_{s} F(X) $ is   through $\hat{\nabla}F(X)(t)=\int_0^{t}\hat{D}_s F(X)\ ds$,  $\forall t\in [0,1]$. Note that $\hat{D}_sF(X)$ is $ds \times d\mu$-almost everywhere well-defined. Then, we have the following relation
$$
    \hat{\nabla}_h F(X)=\langle  \hat{\nabla} F(X), h \rangle_H=\int_0^1 \hat{D}_s F(X)h'_s\ ds= \langle \hat{D} F(X), h' \rangle_{L^2([0,1])}
$$
where $h'$ denotes the derivative of $h$.
\begin{proposition}  \label{chain rule}
Assume $F\in  \M_{p, 1}\left(\mathbb{R}^{d}\right)$, $g: \mathbb{R}^{d} \rightarrow \mathbb{R}$ is a continuous function.
\begin{enumerate}[(i)]
    \item  If $g$ is Lipschitz continuous, then $g\circ F \in \M_{p, 1}$
    \item If $g$ is $\mathcal{C}^1$-function such that
    \begin{align*}
        \E\left[ |g \circ F|^q + \sum_{i}\left|\partial_{i }g\circ F \,  \hat{\nabla} F_i\right|^p \right]<\infty
    \end{align*}
    then $g\circ F \in \M_{r, 1}$ for any $r < p\wedge q$, where $p$ and $q$ are conjugates.
\end{enumerate}
\end{proposition}
\begin{proof}
(i) is evident from Mazur Lemma which says that
closure of a convex set is the same under any topology of the dual pair and from the fact that the graph of $\hat{\nabla}$ is convex in any $L^n(\nu)$, for any $n\geq 1$. \\
(ii) Let $\theta$ be a smooth function of compact support on $\mathbb{R}^d$,  $\theta(0)=1$. Let $\theta_n(x)=\theta(\frac{x}{n})$. Then
$$\E\left[ |\hat{\nabla}(\theta_ng)\circ F|^r\right]\leq 2^{r-1}\sum_i \E\left[|g\circ F|^r|\hat{\nabla} F_i|^r+K |\partial_i g\circ F|^r |\hat{\nabla} F_i|^r\right]
$$
and we have
$$
\E\left[|g\circ F|^r|\hat{\nabla} F_i|^r\right]\leq  \E\Big[|g\circ F|^q \Big]^{r/q}\E\left[ |\hat{\nabla} F_i|^p\right]^{r/p}
$$
where $K$ is an upper bound for $\theta$ and the term with $\theta_n^\prime$ does not contribute. So, $(\theta_n g \circ F,\ n\geq 1)$ is bounded in ${\M}_{r,1}$; hence it has a subsequence which converges weakly and this implies that $\lim_n \theta_n g \circ F=g \circ F$  belongs to ${\M}_{r,1}$.
\end{proof}

\begin{lemma}\label{LemmaFTC}
If $U(X)\in \M_{2,1}(H)$, then we have
\begin{align*}
    \hat{\nabla}_h \int_0^1P(X_s)u_s(X)\cdot dW_s&=\int_0^1 P(X_s)  \hat{\nabla}_h u_s(X)\cdot dW_s + \int_0^1 P(X_s)\partial P(X_s) \hat{\nabla}_h X_s u_s(X)\cdot dW_s \\ &+ \int_0^1 P(X_s)u_s(X)\cdot h'_s \ ds
\end{align*} where  $u_t(X)=U'_t(X)$.
\end{lemma}
\begin{proof}
 Assume that $(u_s)$ is a step process, then
 \begin{align*}
     \int_0^1 P(X_s)u_s(X)\cdot dW_s&=
      \int_0^1 u_s(X)\cdot P(X_s)dW_s  \\
     &= \sum_i u_{s_i}(X)\cdot ( M_{s_{i+1}}- M_{s_{i}})
 \end{align*}
  where $M_t=\int_0^t P(X_s)\ dW_s$ by the martingale representation theorem \cite[Thm.2]{ustunel1}. Therefore, we have
 \begin{align*}
     \hat{\nabla}_h \int_0^1 P(X_s)u_s(X)\cdot  d W_s&=
      \sum_i  \hat{\nabla}_h u_{s_i}(X)\cdot ( M_{s_{i+1}}- M_{s_{i}})\\
     &+ \sum_i u_{s_i}(X)\cdot ({\bf P}(X_{s_{i+1}})h_{s_{i+1}}-{\bf P}(X_{s_i})h_{s_i})\\
     &+\E\left[\sum_i u_{s_i}(X)\cdot \left( \int_{s_{i}} ^{s_{i+1}} \partial P(X_{s})\nabla_h X_{s}\cdot dW_s \right) \mid \mathcal{F}_1(X) \right]
 \end{align*}
 where we define the action of $P(X)$ on the Cameron-Martin space $H$ as ${\bf P}(X_t)h_t=\int_0^t P(X_s)h'_s\ ds$ and use \cite[Prop.2.3.2]{ustunel2010book} for $\nabla M_t$.
 It follows from \cite[Thm. 3]{ustunel2010book} that the last term is equal to
 \begin{align*}
 \sum_i u_{s_i}(X)\cdot \left( \int_{s_{i}} ^{s_{i+1}} P(X_s) \partial P(X_{s})\hat{\nabla}_h X_{s}\cdot dW_s \right)\; .
 \end{align*}
 Then, we pass to the limit in $L^2$ as the mesh of partition goes to zero. For the other terms, the result is straightforward.
\end{proof}

\begin{remark} \label{FTC}
Suppose $u(X_s)$ satisfies the hypothesis of Lemma \ref{LemmaFTC}. Then, we have
$$
\hat{D}_{t} \int_{0}^{1} P(X_s)u(X_s) \cdot d W_s=\int_{0}^{1} P(X_s)  \hat{D}_{t} \left(P(X_s)u(X_s)\right)\cdot d W(s)+P(X_t)u(X_t)
 $$
 $dt\times\mu$-almost everywhere, as $\hat{\nabla}_hF=\int_0^1 \hat{D}_sF h'_s \ ds$.
\end{remark}

\begin{lemma}\label{lemma 2.3}
Let $u(X_s)\in \M_{2,1}(L^2([0,1]))$ be adapted to $\FF(X)$.  Then, $\| u(X_s) \|^2 \in \mathbb{M}_{2,1}$ and
$$
\hat{\nabla}_h \| u(X_s) \|^2 =\hat{\nabla}_h \| U(X) \|_H^2 = 2\langle\hat{\nabla}_h u, U(X) \rangle_H
 $$
 where $U_t(X)=\int_0^t u(X_s)\ ds$.
\end{lemma}
\begin{proof}
The proof is similar to that of Lemma \ref{LemmaFTC}.
\end{proof}
\begin{remark}\label{remark 2.3}
Suppose $u(X_s)$ satisfies the assumption of Lemma \ref{lemma 2.3}. Then, $\int_{0}^{1} |P(X_s)u(X_s)|^2 \ ds$ $ \in \mathbb{M}_{2,1}$ and
$$
\hat{D}_{t} \int_{0}^{1} |u(X_s)|^2 \ ds= 2\int_{0}^{1} \hat{D}_{t}u(X_s)\cdot u(X_s)\ ds.
 $$
\end{remark}

\section{Hedging a Stock Portfolio}

We consider a semimartingale market model with   $n$ risky assets with  price  $X_t=(X_t^1, \ldots,X_t^n)$  and a risk-free asset $X^0_t$ at time $t$.  The asset prices $X_t$ and $X^0_t$    are determined by the system of stochastic differential equations
\begin{align*}
&d X_t=b( X_t )  d t+\sigma( X_t )  d W_{t} \nnumber{model}\\
&d X^{0}_t=r_t X^{0}_t d t
\end{align*}
where $W_t=(W_t^1, \ldots, W_t^n)^*$. We assume that the drift $b$ and the diffusion matrix $\sigma$ satisfy the linear growth and Lipschitz conditions for the existence and uniqueness of a strong solution \cite[Thm. 3.1]{ikeda-watanabe86}. In  \eqref{model}, the arguments of $b$ and $\sigma$ can include time $t$ separately and the analysis of this section will be still valid as this is allowed in our basic reference \cite{ustunel2}, but omitted for the sake of brevity. Examples where the coefficients are functions of only  time $t$ are included in the next section among others. In this section, we will derive the hedging strategy for a given payoff.


Recall that the value process $V^\theta_t$ is given by
\begin{align*}
V^{\theta}_t= \theta^{0}_t\, X^{0}_t +\theta_t\,X_t   \nnumber{value_process}
\end{align*}
where  by $\theta_t, \theta^0_t$ denote the number of shares of $n$ assets and risk-free asset, respectively,  and  $\theta_t$ is taken as a row vector for simplicity of notation. The portfolio $(\theta,\theta^0)$ is  considered to be self-financing, that is,  $V^\theta_t$ satisfies
\begin{align*}
    d V_{t}^{\theta}=\theta^{0}_t d X^{0}_t+\theta_t \,d X_t\nnumber{value_process_diff_0}\; .
\end{align*}
Since from \eqref{value_process}, we have
$
\theta^{0}_t=(V_{t}^{\theta}-\theta_t X_t)/X^{0}_t
$, we rewrite \eqref{value_process_diff_0} as
\begin{align*}
d V_{t}^{\theta}&= r_t\left(V_{t}^{\theta}-\theta_t X_t \right) d t+\theta_t \,d X_t\\
&= \left[ r_t V_{t}^{\theta}+ \theta_t b( X_t ) -r_t  \theta_t X_t\right] d t+\theta_t \sigma(X_t)dW_t \;.
\end{align*}
Assume that the equation below has a solution
\begin{align} \label{marketprice}
    \sigma(X_t) u(X_t)=b( X_t ) -r_t X_t\:.
\end{align}
Although this equation may have several solutions $u$, the orthogonal projection by $P(X_t)$ of these    solutions onto the range space of $\sigma^*(t,X)$ is unique as it can be verified by simple algebra.  Then, the solution $P(X_t) u(X_t)$, called \textit{market price of risk process}, satisfies
\begin{equation}
      \sigma(X_t) P(X_t) u(X_t)=b( X_t ) -r_t X_t\:.            \label{rewrite}
\end{equation}
Note that when \eqref{marketprice} does not have a solution, then the market is not arbitrage-free and this market cannot be used for pricing \cite[pg. 228]{shreve2003}.
Assume also that the $d$-dimensional market price of risk process $u(X_t)$ satisfies
$$
\int_{0}^{1}|P(X_t)u(X_t)|^{2} d t<\infty
$$
almost surely and  the positive local martingale
$$
Z_t \triangleq \exp \left\{-\int_{0}^{t} P(X_s)u(X_s) \cdot d W_s-\frac{1}{2} \int_{0}^{t}|P(X_s)u(X_s)|^{2} d s\right\}
$$
satisfies $\mathbb{E} Z_1 =1$. Then, $Z$ is a martingale with respect to the filtration generated by $X$, $\mathcal{F}(X)$, in view of the converse statement in the martingale representation theorem \cite[Thm.2]{ustunel1}. Now, define $\tilde{W}_t$ by
\begin{align*}
    \widetilde{W}_t=W_t+\int_0^t P(X_s)u(X_s)\ ds.  \nnumber{equl_W}
\end{align*}
Then, $\{\widetilde{W}_t,\ 0\leq t\leq 1 \} $ is  a Brownian motion under the probability measure $\Q$ on $F_1(W)$ given by
\begin{align*}
    \mathbb{Q}(A)=\mathbb{E}[Z_1 1_A],\;\;\; \forall A \in {\cal F}_1(W).
\end{align*}
Using \eqref{rewrite} and\eqref{equl_W}, we can write the price dynamics \eqref{model} using $\widetilde{W}$ as
\[
dX_t= r_t X_t dt + \sigma(X_t)\, d\widetilde{W}_t \; .
\]
Similarly, the discounted price $S_t:= \exp({-\int_0^t r_s ds})X_t$ satisfies
\[
dS_t=e^{-\int_0^t r_s ds}\sigma(X_t)\, d\widetilde{W}_t\; .
\]
Moreover, in view of \eqref{value_process} and \eqref{equl_W}, we can rewrite the value process as
\begin{align*}
d V_{t}^{\theta}= r_t V_{t}^{\theta} d t+\theta_t \sigma(X_t) d \widetilde{W}_t \nnumber{value_process_diff_eq}\;.
\end{align*}
Define the discounted value process
$
     U^\theta_t:= e^{-\int_0^t r_s ds} V^\theta_t
$. Let $G(X)$ be an $\mathcal{F}_1(X)$-measurable and integrable payoff function. After setting  $V_1^\theta=G(X)$ to find the hedging strategy, the equation
\begin{align}  \label{*}
d U^\theta_t=e^{-\int_0^t r_s ds} \theta_t  \sigma(X_t)  \, d \widetilde{W}_t
\end{align}
can be considered as a backward stochastic differential equation with final condition
\begin{equation}
  U_1^\theta =e^{-\int_0^1 r_s ds}G(X)  \;.\label{final}
\end{equation}
Clearly, both the discounted price process and the value process are martingales under $\Q$ when we assume $\int_0^1 \sigma(X_t) \sigma^*(X_t)dt$ is finite a.s. We further assume that
\[
\int \theta_t\sigma(X_t) \sigma^*(X_t)\theta_t^* \, dt < \infty \quad \quad a.s.
\]
to have an admissable strategy $\theta$ \cite[pg.302]{Klebaner}.

\begin{theorem}  \label{thm2}
Suppose $F \in \mathbb{M}_{2,1
}$ is $\mathcal{F}_1(X)$-measurable and the conditions
\begin{align*}
\E [ Z_1^2F^2]<\infty, \quad \quad \E\left[Z_1^2\|\hat{D} F\|^{2} \right]<\infty,
\end{align*}
\begin{align*}
    \E\left[Z_1^2F^2  \left\|\int_{t}^{1}  P(X_s)  \hat{D} \left(P(X_s)u_s\right)\cdot d W(s)+P(X_t)u(X_t) \right. \right.\\
 \left. \left.   +\int_{t}^{1} \hat{D} \left(P(X_s)u(X_s)\right)\cdot P(X_s)u(X_s) d s\right\|^{2}\right]<\infty\;.
\end{align*}
Then, we have  $Z_1F\in \mathbb{M}_{2,1}$ and
\begin{align*}
    F= \E_{\Q}\left[F\right]+\int_0^1   P(X_t)\E_{\Q}\left[\hat{D}_{t} F -F \int_{t}^{1} P(X_s) \hat{D}_{t} \left(P(X_s)u(X_s)\right)\cdot d \widetilde{W}_s \mid \mathcal{F}_{t}(X)\right] d \widetilde{W}_t   \; .
\end{align*}
\end{theorem}
\begin{proof}
We will show that $Z_1F\in \mathbb{M}_{2,1}$ first. Remember that $Z_1=e^{-K}$, where $$K=\int_{0}^{1} P(X_s)u(X_s) \cdot d W_s+\frac{1}{2} \int_{0}^{1} | P(X_s)u (s,X) | ^{2} d s .$$
Lemma \ref{LemmaFTC} implies that $\int_{0}^{1} P(X_s)u(X_s) \cdot d W_s\in \M_{2,1}$ and Lemma \ref{lemma 2.3} implies that $\int_{0}^{1} |P(X_s)u(X_s)|^2 \ ds \in \mathbb{M}_{2,1}$. Hence, $K\in \mathbb{M}_{2,1}.$ Since $\E[F^2e^{2K}],$ $\E[e^{2K}|| \hat{D}F||^2]$ and $\E[F^2e^{2K}|| \hat{D}K||^2]$ are finite by the given assumptions, Proposition \ref{chain rule} implies that $Z_1F\in \mathbb{M}_{2,1}$ satisfying
\begin{align*}
    \hat{D}_t Z_1F&= Z_1 \left( \hat{D}_tF -F (\hat{D}_t K) \right)
\end{align*}
and
\begin{align*}
    \hat{D}_t K&=- \int_{t}^{1} P(X_s)  \hat{D}_{t} \left(P(X_s)u(X_s)\right)\cdot d W_s-P(X_t)u(X_t)  \\
&  \quad \quad -\int_{t}^{1} \hat{D}_{t}(P(X_s)u(X_s))\cdot P(X_s)u_s ds\; .
\end{align*}
Let
$
Y_t= \E_{\mathbb{Q}}\left[F \mid \mathcal{F}_{t}(X)\right]
$ and note that
\begin{align*}
Z^{-1}_t &=\exp \left\{\int_{0}^{t} P(X_s)u(X_s) \cdot d W_s+\frac{1}{2} \int_{0}^{t} | P(X_s)u (s,X) | ^{2} d s\right\}\\
&=\exp \left\{\int_{0}^{t}  P(X_s)u(X_s) \cdot  d \widetilde{W}_s-\frac{1}{2} \int_{0}^{t} | P(X_s)u (s,X) | ^{2} d s\right\}\nnumber{eq_Z_tilde}
\end{align*}
 Then, we get
\begin{align*}
Y_{t} &=Z^{-1}_t \E\left[Z_1 F \mid \mathcal{F}_{t}(X)\right]
\\&=Z^{-1}_t\left\{\E[Z_1 F]+\int_{0}^{1} P(X_s)\E\left[\hat{D}_{s}
\E\left[Z_1 F \mid \mathcal{F}_{t}(X)\right] \mid \mathcal{F}_{s}(X)\right] \cdot d W_s\right\}
\\
&=Z^{-1}_t\left\{\E[Z_1 F]+\int_{0}^{t} P(X_s)\E\left[\hat{D}_{s}(Z_1 F) \mid \mathcal{F}_{s}(X)\right] \cdot d W_s\right\}  =:Z^{-1}_t A_t\nnumber{eq_Y}
\end{align*}
where we have applied the formula of \cite[Thm.6]{ustunel2}  to $\E\left[Z_1 F \mid \mathcal{F}_{t}(X)\right]$ and used the fact $\hat{D}_{s}
\E\left[Z_1 F \mid \mathcal{F}_{t}(X)\right]$ is $\FF_t(X)$-measurable for $t>s$ and equal to $0$ otherwise.
From \eqref{eq_Z_tilde} and \eqref{eq_Y}, we get
\begin{align*}
    dZ^{-1}_t& = Z^{-1}_t P(X_t)u(X_t) \ d\widetilde{W}_t\\
    dA_t &= P(X_t)\E\left[\hat{D}_{t}(Z_1 F) \mid \mathcal{F}_{t}\right] \cdot d W_t\\
    dA_tdZ^{-1}_t&= Z^{-1}_t P(X_t)u(X_t) P(X_t)\E\left[\hat{D}_{t}(Z_1 F) \mid \mathcal{F}_{t}\right] \ dt
\end{align*}
Since $
    d Y_t=A_t\ dZ^{-1}_t +Z^{-1}_t\ dA_t + dA_t\ dZ_t^{-1}
$, it follows that
\begin{align*}
d Y_t&= \left\{\E[Z_1 F]+\int_{0}^{t} P(X_s)\E\left[\hat{D}_{s}(Z_1 F) \mid \mathcal{F}_{s}(X)\right] \cdot d W_s\right\}Z^{-1}_t P(X_t)u(X_t) \cdot d\widetilde{W}_t\\
&\quad + Z^{-1}_t P(X_t)\E\left[\hat{D}_{t}(Z_1 F) \mid \mathcal{F}_{t}(X)\right]\cdot d W_t\\
&\quad+Z^{-1}_t P(X_t)u_t  P(X_t)\E\left[\hat{D}_{t}(Z_1 F) \mid \mathcal{F}_{t}(X)\right]\ dt\\
&= Y_t P(X_t)u_t d\tilde{W}_t  +  Z^{-1}_t P(X_t)\E\left[\hat{D}_{t}(Z_1 F) \mid \mathcal{F}_{t}(X)\right] d \widetilde{W}_t
\\
&=P(X_t)u_t \E_{\Q}\left[F \mid \mathcal{F}_{t}(X)\right] d \widetilde{W}_t
\\
&+P(X_t)\E_{\Q}\left[\hat{D}_{t} F \mid \mathcal{F}_{t}(X)\right]\widetilde{W}_t
\\
&+P(X_t)\E_{\Q}\left[F \left\{ -\int_{t}^{1} P(X_s)  \hat{D}_{t} \left(P(X_s)u_s\right)\cdot d W(s)- P(X_t)u_t  \right\} \mid \mathcal{F}_{t}(X)\right] d \widetilde{W}_t
\\
&+P(X_t)\E_{\Q}\left[F \left\{   -\int_{t}^{1} \hat{D}_{t} (P(X_s)u(X_s))\cdot P(X_s)u(X_s)  ds \right\} \mid \mathcal{F}_{t}(X)\right] d \widetilde{W}_t
\\
&=P(X_t)\E_{\Q}\left[\hat{D}_{t} F -F \int_{t}^{1} P(X_s) \hat{D}_{t} \left(P(X_s)u(X_s)\right)\cdot d \widetilde{W}_s \mid \mathcal{F}_{t}(X)\right] d \widetilde{W}_t
\end{align*}
In view of  $
Y_1=    \E_{\Q}\left[F \mid \mathcal{F}_{1}\right]=F $ and
$Y_0= \E_{\Q}\left[F \mid \mathcal{F}_{0}\right]= \E_{\Q}\left[F\right]
$,
we get
\begin{align*}
    F= \E_{\Q}\left[F\right]+\int_0^1   P(X_t)\E_{\Q}\left[\hat{D}_{t} F -F \int_{t}^{1} P(X_s) \hat{D}_{t} \left(P(X_s)u(X_s)\right)\cdot d \widetilde{W}_s \mid \mathcal{F}_{t}(X)\right] d \widetilde{W}_t   \; .
\end{align*}
\end{proof}

Now, we are ready to find the hedging strategy for our market model, when a payoff function $G$ is given. Letting  $F:=U_1^\theta  =e^{-\int_0^1 r_s ds}G(X)$ in \eqref{final}, which needs to hold for the aim of finding a replicating portfolio, and substituting $F$ in the result of Theorem \ref{thm2} with the assumption that its conditions are satisfied, we get
\begin{align*}
&  U_1^\theta =  \; \E_\Q[e^{-\int_0^1 r_s ds}G]+ \nnumber{dis_value_process_clark} \\
&  \int_0^1   P(X_t)\E_{\Q}\left[\hat{D}_{t} (e^{-\int_0^1 r_s ds}G) -e^{-\int_0^1 r_s ds}G \int_{t}^{1} P(X_s) \hat{D}_{t} \left(P(X_s)u(X_s)\right)\cdot d \widetilde{W}_s \mid \mathcal{F}_{t}(X)\right] d \widetilde{W}_t   \; .
\end{align*}
On the other hand, in view of \eqref{*} and as $U_0^\theta=V_0^\theta$ by definition, we have
\begin{align*}
     U^\theta_1&=V^\theta_0+ \int_0^1  e^{-\int_0^t r_s ds} \theta_t \sigma(X_t)   \, d \widetilde{W}_t  \; .
\end{align*}
In comparison with \eqref{dis_value_process_clark}, we conclude   that
\begin{align}
    V_0^\theta= \E_\Q[e^{-\int_0^1 r_s ds}G(X)] \label{price}
\end{align}
and the hedging strategy $\theta_t$  solves
\begin{align*}\nnumber{Hedging_str_long}
   & \sigma^*(X_t) \theta_t^* =e^{\int_0^t r_s ds} P(X_t)\\
    & \quad \quad  \E_{\Q}\left[\hat{D}_{t} (e^{-\int_0^1 r_s ds}G(X)) -e^{-\int_0^1 r_s ds}G(X) \int_{t}^{1} P(X_s) \hat{D}_{t} \left(P(X_s)u(X_s)\right)\cdot d \widetilde{W}_s \mid \mathcal{F}_{t}(X)\right]
\end{align*}
at each time $t\ge 0$, where $r_s$ can be understood as a function of the asset prices and denoted as $r(X_s)$ if it is random.
Note that $u$ appears in $\Q$ and   $\widetilde{W}$ in the above equation, and we can obtain a unique and adapted strategy if we replace $\theta_t$ by $P(X_t)\theta_t$.
If the interest rate is non-random, then Equation \eqref{Hedging_str_long}  reduces to
\begin{align*}\nnumber{Hedging_str}
   & \sigma^*(X_t) \theta_t^* =e^{-\int_t^1 r_s ds} P(X_t)\\
    & \quad \quad \quad \quad \E_{\Q}\left[\hat{D}_{t} G(X) -G(X) \int_{t}^{1} P(X_s) \hat{D}_{t} \left(P(X_s)u(X_s)\right)\cdot d \widetilde{W}_s \mid \mathcal{F}_{t}(X)\right].
\end{align*}

\begin{remark}
 If $\sigma$ is non-degenerate, then the projection map $P(X_s)$ is just the identity map, $\FF(X)=\FF(W)$ and $\hat{\nabla}=\nabla$. Assuming that $\sigma$ is a square matrix for simplicity, we can rewrite  \eqref{Hedging_str} as
\begin{align*}
     \theta_t^* =\sigma^*(X_t)^{-1}\ e^{\int_0^t r_s ds}\ \E_{\Q}\left[{D}_{t}(e^{-\int_0^1 r_s ds}G(X)) - e^{-\int_0^1 r_s ds}G(X) \int_{t}^{1}  {D}_{t} \left(u(X_s)\right)\cdot d \widetilde{W}_s \mid \mathcal{F}_{t}(W)\right]
\end{align*}
which is the same as  \cite[Eq.(3.10)]{ocone1991}.
\end{remark}

\section{Examples for Payoff Function}

In this section, the hedging strategy is worked out for some examples of the payoff function to demonstrate the formulas. Here, we indicate  the  terminal time by $T$. From the point of view of option pricing with a claim $G(X)$, the analysis of a hedging strategy $\theta$ can  readily be used. The claim $G(X)$ is attainable if $\mathbb{E} [G(X)] < \infty$ and there exists an admissable strategy $\theta_t$, $0\le t\le 1$ \cite[pg.303]{Klebaner}. Then, the price of the claim at time $t$ is given by
 \[
 e^{\int_t^1 r_s ds}\,\E_\Q[G\mid \mathcal{F}_{t}(X)]
 \]
with the assumption that the interest rate $r$ is deterministic, and in particular at time 0, the price is equal to \eqref{price}. We consider various claims below as suitable for demonstration of our results.

\subsection{Black–Scholes model}
Consider the one-dimensional Black–Scholes model
\begin{align*}
dX ^{0}_t&=\rho X^{0}_t d t,\ X^{0}_0=1\\
d X^1_t&=\mu X^{1}_t d t+\sigma X^{1}_t d W^1_t,\ X_{0}^1>0
\end{align*}
where $\rho,\mu, \sigma >0$. The equivalent martingale measure for this one-dimensional model is  $
\mathbb{R}(A)=\E[Y_1 1_A]
$,
where
\begin{align*}
    Y_t=\exp \left\{-\int_{0}^{t} u\ d W_s^1-\frac{1}{2} \int_{0}^{t} u^{2}\ d s\right\}, \quad 0 \leq t \leq 1
\end{align*}
and $u =(\mu-\rho)/\sigma$.
 For this model, the hedging strategy is given by
\begin{align}  \label{oks}
   \theta_{t}=e^{\rho t} \sigma^{-1} (X^1_t)^{-1} \E_{\mathbb{R}}\left[D_{t} G \mid \mathcal{F}_{t}(W^1)\right]
\end{align}
where $D_t G=\frac{\text{d}}{dt}\nabla G$ and $\nabla$ is Gross-Sobolev derivative defined for the functionals of $W^1$ \cite[Ex.4.1.1]{nunno_2009}.
Clearly, this model is not an example of the degenerate case. However, we can rewrite it as a degenerate model by artificially taking $X^2:=X^0$ to demonstrate our formulas. In \eqref{model}, take   $ r_t=\rho $
\begin{align*}
    b(X_t)=(\mu X^1_t,\rho X_t^2)
\quad \quad  \sigma(X_t)=\begin{bmatrix}
\sigma X^1_t & 0 \\
0 & 0
\end{bmatrix}.
\end{align*}
Observe that
\begin{align*}
    P(X_t)=\begin{bmatrix}
 1 & 0 \\
0 & 0
\end{bmatrix} \quad \quad
P(X_t)u(X_t)=\begin{bmatrix}
 \frac{\mu-\rho}{\sigma} \\
0
\end{bmatrix} \quad \quad
\hat{D}_t P(X_t)u(X_t)=(0,0).
\end{align*}
If we substitute these in  \eqref{Hedging_str}, we get
\begin{align*}
    \begin{bmatrix}
\sigma X^1_t & 0 \\
0 & 0
\end{bmatrix}\theta_t^* =e^{\rho t} \begin{bmatrix}
 1 & 0 \\
0 & 0
\end{bmatrix}\E_{\Q}\left[\hat{D}_{t} G  \mid \mathcal{F}_{t}(X)\right].
\end{align*}
Moreover, it is easy to see that
$\FF_t(X)=\FF_t(W^1),$  $\Q(A)=\mathbb{R}(A),$ for each $A\in \FF_t(X)$, and $\hat{D}_tG=(D_t G,0),$ where derivative $D_t G$ taken in the sense of Malliavin calculus for Brownian motion. Hence, the hedging strategy is $\theta_t^*=[\theta_t,0] $ with $\theta_t$ of \eqref{oks}. When $G$ is taken to be the European option $G=(X_T^1-K)^+$, we have $G\in \M_{2,1}$ by Proposition \ref{chain rule}.
The hedging portfolio for $G$ is given by \begin{align*}
\theta_t^* = \begin{bmatrix}
 e^{\rho t}\sigma^{-1}(X^1_t)^{-1}\E_{\mathbb{R}}[D_{t} F  \mid \mathcal{F}_{t}(X)]  \\
0
\end{bmatrix}
\end{align*}
in this case, equivalent to the result in \cite[Ex 4.1.1]{nunno_2009}.
\subsection{Exotic  Options}
Exotic options are a class of options contracts in that the value of an option and the time that the holder can exercise it  depend on the prices of the assets on the whole period  \cite{Klebaner}. Since exotic options can be  customized to the needs of the investor, it provides various investment alternatives. We will examine exotic options in a two-dimensional market model with terminal time $T$. Without loss of generality, we assume that $\sigma_{11}$ in \eqref{model} is away from zero. Let $\Delta=\sigma_{11}(X_s)\sigma_{22}(X_s)-\sigma_{12}(X_s)\sigma_{21}(X_s)$.
	When $\Delta=0$, the projection map  can be written as
	$$
	P\left(X_s\right)=\frac{1}{\sigma_{11}^2(X_s)+\sigma_{12}^2(X_s)}\left[\begin{array}{cc}
	\sigma_{11}^2(X_s) & \sigma_{11}(X_s) \sigma_{1 2}(X_s) \\
	\sigma_{11}(X_s) \sigma_{12}(X_s) & \sigma_{12}^2(X_s)
	\end{array}\right]
	$$
	and  the projected market price of risk process is given by
	$$
	P\left(X_{s}\right) u\left(X_{s}\right)=\frac{b_{1}\left(X_{s}\right)-r_{s} X_{s}^{1}}{\sigma_{11}^{2}\left(X_{s}\right)+\sigma_{1 2}^{2}\left(X_{s}\right)}\begin{bmatrix}
	\sigma_{11}\left(X_{s}\right) \\
	\sigma_{1{2}}\left(X_{s}\right)
	\end{bmatrix}
	$$
	in view of \eqref{marketprice}.
	Suppose $b$ and $\sigma$ have bounded partial derivatives, then $P\left(X_{s}\right) u\left(X_{s}\right)$ has the derivative
	$\hat{D}_t P\left(X_{s}\right) u\left(X_{s}\right)=J_f(X_s)(\hat{D}_t X^1_s,\hat{D}_t X^2_s)$, where $J_f$ is the Jacobian of $f(x,y)=\frac{b_{1}\left(x,y\right)-r_{s} x}{\sigma_{11}^{2}\left(x,y\right)+\sigma_{1 2}^{2}\left(x,y\right)}\begin{bmatrix}
	\sigma_{11}\left(x,y\right) \\
	\sigma_{1{2}}\left(x,y\right)
	\end{bmatrix}$ and $\hat{D}_t X_s$ solves
	\begin{align*}
	\hat{D}_t X_s = \int_t^s  J_b(X_r) \hat{D}_t(X_r)\ d r+\int_t^s   J_{\sigma_i } \hat{D}_t(X_r)\cdot P(X_r)d W_{r} + \sigma(t,X_t)
	\end{align*}
	$dt \times d\mu$ - a.e., where $\sigma_i$ is $i^{th}$ row of matrix $\sigma(X_s)$ \cite[Thm. 5]{ustunel2}, and $J_b$, $J_{\sigma_i}$ denote  the Jacobian matrices of $b$ and $\sigma_i$, respectively.
	
	Consider the linear case, that is, $b(X_s)=(b_1(s)X^1_s, \ b_2(s)X^2_s)$ and
	\begin{align*}
	\sigma =\left[\begin{array}{cc}
	\sigma_{11}(s)X^1_s & \sigma_{12}(s)X^1_s \\
	\sigma_{21}(s)X^2_s & \sigma_{22}(s)X^2_s
	\end{array}\right]  \; .
	\end{align*}
	Then, we have
	$$
	P\left(X_{s}\right) u\left(X_{s}\right)=\frac{b_{1}\left({s}\right)-r_{s} }{\sigma_{11}^{2}\left({s}\right)+\sigma_{12}^{2}\left({s}\right)}\left[\begin{array}{c}
	\sigma_{11}\left({s}\right) \\
	\sigma_{1{2}}\left({s}\right)
	\end{array}\right].
	$$
	Clearly, $P\left(X_{s}\right) u\left(X_{s}\right)$ is deterministic, which implies
	$\hat{D}_tP(X_{s}) u(X_{s})=0$ and the hedging strategy $\theta_t$ at time $t$ solves
	$$\sigma^*(X_t) \theta_t^* =e^{-\int_t^T r_s ds} P(X_t)\E_{\Q}\left[\hat{D}_{t} G \mid \mathcal{F}_{t}(X)\right]$$ for given payoff function $G$ by \eqref{Hedging_str}.
	
{\em 	Asian options} are options where the payoff depends on the average of the underlying
	assets. Pricing of Asian options have been studied in \cite{Rogers95,vecer2002} by the use of PDEs. In \cite{Yang2011}, Yang et al. have used Malliavin calculus to derive the
	hedging strategy and price of Asian option. We will consider Asian call option with floating strike, which
	  pays at time $T,$ the payoff $G(X)= \left(\frac{1}{T} \int_0^T X_s^1\ ds- K X_T^2 \right)^+$. Proposition \ref{chain rule} implies $G\in\M_{2,1}$ and  $\hat{D}_t G(X)=1_A \left(\frac{1}{T} \int_t^T \hat{D}_t X_s^1\ ds -K\hat{D}_t X_T^2 \right)$ for  $dt \times d\mu$ - a.e., where $A=\{ \frac{1}{T} \int_0^T X_s^1\ ds- K X_T^2>0 \}$
	and
	$$
	\hat{D}_{t} X_{s}^i=\begin{bmatrix}
	\sigma_{i1}\left(s\right) \\
	\sigma_{i2}\left(s\right)
	\end{bmatrix}\exp \left(\int_{t}^{s}P(X_r)\sigma_i(r) \cdot d W_{r}+\int_t^{s} b_i(r)\ d r-\frac{1}{2} \int_{t}^{s} |P(X_r) \sigma_i(r)|^2\ d r  \right).$$
	
	{\em An exchange option}  gives the  right to put a predefined risky
asset and call the other risky asset, as introduced in \cite{Margrabe1978}. The price and hedging strategy have been calculated in \cite{Mataramvura2012} via Malliavin calculus. The  payoff function $G(X)=(X^1_T-X^2_T)^+$ has the derivative given by
	$\hat{D}_t G(X)=1_A \left( \hat{D}_t X_T^1 -\hat{D}_t X_T^2 \right)$ for $dt \times d\mu$ - a.e., where $A=\{  X_T^1- X_T^2>0 \}$.
	
	{\em Look-back options} are a particular type of path-dependent options so that the value of the payoff function depends on the minimum or maximum of the underlying asset price.
The price of option and the hedging strategy have been derived in \cite{Bermin2000,Renaud2007} with Malliavin calculus. Defining
$M_{0, T}^{X^{1}}=\sup_{0\le t \le T}X_t^1 $, we  consider the payoff  $G(X)=\left(M_{0, T}^{X^{1}} -KX^2_T\right)^+$, which has the derivative
 \begin{align*}
 \hat{D}_{t} G(X)=& 1_A \left(\hat{D}_{t}M_{0, T}^{X^{1}}  -K\hat{D}_t X_T^2\right)
 \end{align*}
 with
 \begin{align*}
    \hat{D}_{t}M_{0, T}^{X^{1}} =\begin{bmatrix}
    \sigma_{i1}(\tau)\\ \sigma_{i2}(\tau)
    \end{bmatrix}\exp
    \left(\int_{t}^{\tau}P(X_s)\sigma_i(s) \cdot d W_{s}+\int_t^{\tau} b_i(s)\ d s-\frac{1}{2} \int_{t}^{\tau} |P(X_s) \sigma_i(s)|^2\ d s  \right)
\end{align*}
for 	$dt \times d\mu$ - a.e.,  where
 $\tau=\inf \{t:\ X_t^1=M_{0, T}^{X^{1}} \}$ and $A=\{M_{0, T}^{X^{1}} -KX^2_T>0\}$. Here, we have used the approach in \cite[pg. 55]{Li_PhD2011} to calculate the derivative of $M_{0, T}^{X^{1}}$ as follows. The model considered in \cite{Li_PhD2011}  is one-dimensional, has constant volatility and Malliavin calculus for Brownian motion is used to calculate the hedging portfolio. However, the idea can be applied easily to our case. For each $m \in \mathbb{N},$ choose a partition $\pi_m=\left\{0=s_{1}<\cdots<s_m=T\right\}$ so that $\pi_m \subseteq \pi_{m+1}$ and $\cup_m \pi_m$ is dense in $\left[0, T\right]$. Define $\varphi_m$ by $\varphi_m(\underline{x})=\max _{1 \leq i \leq m} x_{i}$ . Then, $\varphi_m(X^1_{\pi_m})\rightarrow M_{0,T}^{X^1}$ in $L^2$. Since the function $\varphi_m$ is Lipschitz, $\varphi_m(X^1_{\pi_m})\in \M_{2,1}$ by \ref{remark 2.3}. Moreover, $\sup_m\| \hat{D}_t \varphi_m(X^1_{\pi_m}) \|^2 \leq \sup_t \sigma_t T \|M_{0,T}^{X^1}\|^2_{L^2(\mu)}$ and
this implies $M_{0,T}^{X^1}\in \M_{p,1}$. Note that
$$\hat{D}_t \varphi_m(X^1_{\pi_m}) = \begin{bmatrix}
    \sigma_{i1}(\tau_m)\\ \sigma_{i2}(\tau_m)
    \end{bmatrix}\exp
    \left(\int_{t}^{\tau_m}P(X_s)\sigma_i(s) \cdot d W_{s}+\int_t^{\tau_m} b_i(s)\ d s-\frac{1}{2} \int_{t}^{\tau_m} |P(X_s) \sigma_i(s)|^2\ d s  \right) $$
    where ${\tau_m}=\min \left\{t_{j}: X^1_{t_{j}}=\varphi_m(X^1_{\pi_m})\right\}$. For each $m,$ $\tau_m$ is a measurable function and $\tau_m\rightarrow\tau$ a.s. When  $\sigma_{i1}$ is assumed to be continuous, the result follows as      $\hat{D}$ is a closed operator.

	
\section{Conclusion}

We have used Malliavin calculus for degenerate diffusions to derive the hedging portfolio for a given payoff function in a semimartingale market model. The prices are assumed to follow a multidimensional diffusion process with  a singular volatility matrix $\sigma \sigma^*$, where $\sigma$ is taken to be a function of the prices  with no extra randomness. In applications, the estimation of the volatility $\sigma \sigma^* $  is crucial from financial data, which may be accomplished for example through estimation of the variance of the price time series, yielding an estimate for $\sigma \sigma^*$. In the case of degeneracy of the estimate, this can now be taken care of by the results of the present paper.

From a theoretical point of view, the projection operator $P$ to the range space of  $\sigma^*$ plays a crucial role in our results. We have shown that the hedging portfolio can be obtained uniquely as a solution to a system of linear equations by projecting any solution of the system with $P$. For this purpose, a version of the Clark-Ocone type formula for functionals of degenerate diffusions is proved  under an equivalent change of measure.  As  demonstration  of our results,   intermediate  calculations for the Gross-Sobolev type derivative of the payoff function of the  prices are given in the case of exotic  options in  a two-dimensional linear model.

\bibliographystyle{plain}

\bibliography{Hedging.bib}
\end{document}